\documentclass[12pt]{article}
\usepackage{authblk}
\usepackage{amsfonts, amsmath, amssymb, amsthm, mathtools, relsize}
\usepackage{hyperref}
\usepackage{graphicx}
\usepackage[normalem]{ulem}
\usepackage{fullpage}
\usepackage{enumerate}
\usepackage{graphicx, float}
\usepackage{color}

\usepackage[numbers,sort&compress]{natbib}

\usepackage{tikz,capt-of}
\usetikzlibrary{automata, positioning}
\usetikzlibrary{shapes.geometric, arrows, shapes.misc}

\title{Distributions of cherries and pitchforks} 

\newcommand{\al}{\alpha}

\newcommand{\tw}[1]{{\textcolor{black}{#1}}}

\newcommand{\FF}{\mbox{${\mathcal F}$}}
\newcommand{\NN}{\mbox{${\mathcal N}$}}

\newcommand{\OO}{\mbox{${\mathcal O}$}}

\newcommand{\Ebold}{\mbox{${\mathbb E}$}}

\newcommand{\pp}{\mathbb{P}}

\newcommand{\bu}{{\mathbf u}}
\newcommand{\bv}{{\mathbf v}}

\newcommand{\var}{\mathrm{var}}
\newcommand{\cov}{\mathrm{cov}}

\newcommand{\cas}{\ \stackrel{\mbox{a.s.}}{\longrightarrow} \ }

\newcommand{\Ga}{\Gamma}       
\newcommand{\si}{\sigma}       

\newcommand{\pf}{A}
\newcommand{\ch}{C}
\newcommand{\tU}{\widetilde{U}}

\newcommand{\udim}{d}

\newcommand{\proj}{T}
\newcommand{\conv}{Q}

\newcommand{\nal}{\beta}
\newcommand{\wconv}{\xrightarrow{~d~} }
\newcommand{\bzero}{\mathbf{0}}

\newcommand{\pev}{\lambda_1} 

\newcommand{\gp}{\phi}  

 \newtheorem{theorem}{Theorem}[section]
 \newtheorem{proposition}[theorem]{Proposition}
 \newtheorem{corollary}[theorem]{Corollary}
 
 \newtheorem{lemma}[theorem]{Lemma}

 \newtheorem{remark}{Remark}

\begin{document}
\title{Distributions of cherries and pitchforks for the Ford model}

\author{Gursharn Kaur\thanks{Biocomplexity Institute, University of Virginia, Charlottesville, USA 22911.}}
\author{Kwok Pui Choi\thanks{Department of Statistics and Data Science,  and the Department of Mathematics, National University of Singapore, Singapore 117546.}}
\author{Taoyang Wu\thanks{School of Computing Sciences, University of East Anglia, Norwich, NR4 7TJ, U.K.}}

\affil[]{}

	 \maketitle

\begin{abstract}
We study two fringe subtree counting statistics, the number of cherries and that of pitchforks for  Ford's $\alpha$ model, a one-parameter family of random phylogenetic tree models that includes the uniform and the Yule models, two tree models commonly used in phylogenetics. Based on a nonuniform version of the extended P\'olya urn models in which negative entries are permitted for their replacement matrices, we obtain the strong law of large numbers and the central limit theorem for the joint distribution of these two count statistics for the Ford model.  Furthermore, we derive a recursive formula  for computing the exact joint distribution of these two statistics. This leads to exact formulas for their means and higher order asymptotic expansions of their second moments, which allows us to identify a critical parameter value for the correlation between these two statistics. That is, when $n$ is sufficiently large,  they are negatively correlated for $0\le \alpha \le 1/2$ and positively correlated for $1/2<\alpha<1$.
\end{abstract}

\section{Introduction}
In biology,  evolutionary relationships among the biological system under investigation are typically represented by a phylogenetic tree, that is, a binary tree whose leaves are labelled by a set of species. 
Distributional properties of tree shape statistics, such as Colless' index, Sackin's index and the number of subtrees, under random tree models play an important role in investigating phylogenetic trees inferred from real datasets~(see, e.g.~\cite{steel2016phylogeny}). Various properties concerning  these  statistics have been established in the past decades on the following two fundamental random phylogenetic tree models: 
the Yule model (aka the Yule-Harding-Kingman (YHK) model)~\cite{rosenberg06a,disanto2013exact,Janson2014} and the uniform model (aka the proportional to distinguishable arrangements (PDA) model)~\cite{McKenzie2000,chang2010limit, WuChoi16,CTW19}.
However, for phylogenetic trees inferred from real datasets, the Yule or uniform model may not always be a good fit~\cite{blum2006random}, and several general classes of random trees have been proposed  for modelling and analysing the observed data,
two popular ones being Ford's alpha model~\cite{Ford2006} and Aldous' beta model~\cite{aldous96a}.

In this paper,  we confine ourselves to  Ford's alpha model, a one-parameter family  of random tree growth models introduced by Daniel J. Ford in his PhD thesis \cite{Ford2006}. 
 More precisely, under the Ford model with a fixed parameter $0\le \alpha \le 1$, a random tree of a given number of leaves is generated such that at any step in which a tree $T_n$ with $n$ leaves has been constructed from previous steps, a new leaf attaches to  an internal edge of $T_n$ with probability $\frac{\al}{n-\al}$ and  to  a leaf edge in $T_n$ with probability $\frac{1-\al}{n-\al}$. The resulting random tree model will be referred to as the Ford model (indexed by the parameter $\alpha$) in this paper, which is also known as the alpha tree model~(see, e.g.~\cite{coronado2019balance}). Note that the Ford model is a family of random tree models which includes the Yule model with $\al=0$ (which is closely related to random binary search trees, see e.g.~\cite{Janson2014}),   the uniform model with $\al=1/2$, and the Comb model with $\al=1$.  

The tree shape indices studied in this paper are the number of cherries and the number  of pitchforks. A cherry is a fringe subtree (i.e., a subtree consisting of one of the edges $(u,v)$ and all the descendants of $v$) with precisely two leaves and a pitchfork a fringe subtree with three leaves. 
The study of the number of fringe subtrees of a random tree can be traced back to a paper of Aldous~\cite{aldous1991asymptotic} and has since been extended to various random tree models
(see, e.g.~\cite{holmgren2017fringe}).
In phylogenetics, the asymptotic properties of the number of cherries was first studied by McKenzie and Steel  \cite{McKenzie2000}, who showed that the number of cherries is asymptotically normal for the Yule and the uniform models as the number of leaves tends to infinity. Later, similar properties of the number of cherries are  extended to the Ford model~\cite[Theorem 57]{Ford2006} and to the Crump-Mode-Jagers branching process~\cite{plazzotta2016asymptotic}. 
For the number of pitchforks, Rosenberg~\cite{rosenberg06a} obtained its  mean and variance, and Chang and Fuchs~\cite{chang2010limit} proved that  the number of pitchforks is also  asymptotically normal for the Yule and the uniform models. For the joint distributions, Holmgren and Janson showed that~\cite{Janson2014} the joint distribution is asymptotically normal for the Yule model, using a correspondence between the Yule model and random binary search trees. This was recently extended to the uniform model based on a uniform version of the extended urn models in which negative entries are permitted for their replacement matrices~\cite{Paper1}.

In this paper,  we establish the strong law of large numbers and the central limit theorem 
for the joint distribution of cherries and pitchforks under the Ford model (Theorem~\ref{Thm:Convg-ChPh}) 
by considering an associated nonuniform urn model (Theorem~\ref{thm:urn:edge}). These results are presented in Section~\ref{sec:limiting}, following Section~\ref{sec:preliminary} in which we  collect  background concerning the Ford model and  limiting theorems on uniform urn models.

\tw{
Furthermore, we derive a recurrence  formula for computing the exact joint distribution under the Ford model~(Theorem~\ref{jointpmf}) in Section~\ref{sec:exact}, generalizing the results in~\cite{WuChoi16} for the Yule and the uniform models. This enables us to obtain recurrence expressions for 
the moments of these two statistics. 
As an application, in Theorems~\ref{thm:firstm} we obtain the exact formula for the mean of the number of cherries (first reported in~\cite{Ford2006}) and that of  pitchforks under the Ford model. Furthermore, we also obtain 
higher order expansions of the second moments of their marginal and joint distributions~(Theorems~\ref{thm:secondm}), which allows us to identify a critical parameter value for the correlation between these two statistics. That is, when $n$ is sufficiently large these two statistics are negatively correlated for $0\le \alpha \le 1/2$ and positively correlated for $1/2<\alpha<1$. 
The proofs of Theorems~\ref{thm:firstm}~\&~\ref{thm:secondm} are presented in   Section~\ref{sec:expansion}.
}


\section{Ford Model and Urn Model}
\label{sec:preliminary}
In this section, we first introduce  the Ford model, which is a one-parameter family of random phylogenetic tree models. Next we present a nonuniform version of the extended urn models  associated with the Ford tree model. Finally, we recall certain conditions on the related uniform version of the extended urn model under which the strong law of large numbers and the central limit theorem are obtained.

\subsection{Ford model}

A rooted binary tree is a finite connected simple graph without cycles that contains a unique vertex of degree 1 designated as the root and all the remaining vertices are of degree 3 (interior vertices) or 1 (leaves).  A phylogenetic tree  $T$ with $n$ leaves is a rooted binary tree whose leaves are bijectively labelled by the elements in $\{1,\dots,n\}$. Note that all the edges in $T$ are directed away from the root and edges incident with leaves are referred to as pendant edges. A fringe subtree in $T$ consisting of one of the edges $(u,v)$ and all the descendants of $v$. A cherry (resp. pitchfork) is a fringe subtree with two (resp. three) leaves. A cherry not contained in a pitchfork will be referred to as an essential cherry. Finally, we let $A(T)$ and $C(T)$ denote the number pitchforks and that of cherries in tree $T$, respectively.

Under the Ford model with parameter $0\le \alpha \le 1$,  a random phylogenetic tree $T_n$ with $n$ leaves is constructed recursively by adding one leaf at a  time as follows. Fix a random permutation $(x_1,\dots,x_n)$ of $\{1,\dots,n\}$.  The initial tree $T_2$  contains precisely two leaves (e.g. one cherry) which are labelled as $x_1$ and $x_2$. For the recursive step, given a tree $T_m$ with $m$ leaves constructed so far, choose a random edge in $T_m$ according to the distribution that assigns weight $1-\al$ to each pendant edge (i.e., those incident with a leaf)  and weight $\al$ to each of the other edges. The new leaf labelled $x_{m+1}$ bifurcates the selected edge and joins in the middle. 
Every single addition of a leaf in the tree results into a replacement of the selected edge with two new edges.
Finally, we let $A_n=A(T_n)$ and $C_n=C(T_n)$ denote the numbers of pitchforks and cherries in tree $T_n$, respectively.

 \begin{figure}[ht]
	\begin{center}
		{\includegraphics[width=0.9\textwidth]{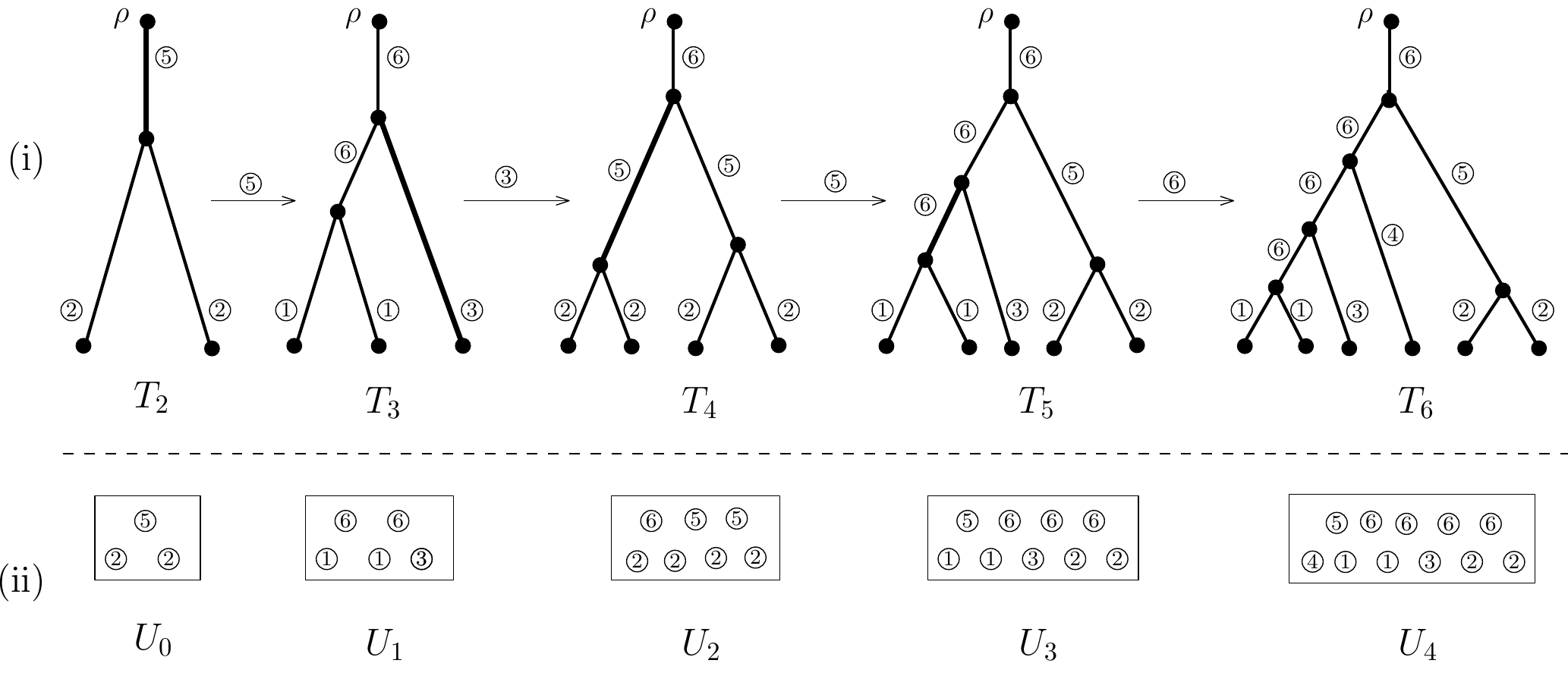}}
	\end{center}
	\caption{A sample path of the Ford model and the associated trajectory under the urn model. (i) A sample path of the Ford model evolving from $T_2$ with two leaves to $T_6$ with six leaves. The labels of the leaves are omitted for simplicity. The type of an edge is indicated by the circled number next to it. For $2\le i \le 5$, the edge selected in $T_i$ to generate $T_{i+1}$ is highlighted in bold and the associated edge type is indicated in the circled number above the arrow. (ii) The associated urn model with six colours, derived from the types of pendants edges in the trees. 
		In vector form, $U_0=(0,2,0,0,1,0), U_1=(2,0,1,0,0,2), U_2=(0,4,,0,0,2,1), U_3=(2,2,1,0,1,3)$, and $U_4=(2,2,1,1,1,4)$. 	
	}
	\label{fig:example}
\end{figure}

\subsection{An urn model associated with trees }
\label{subsect:urn:tree}

Consider an urn containing balls of $d$ different colours where colours are denoted by integers $\{1,2,\dots, \udim\}$. Let $U_n=(U_{n,1},\dots, U_{n,\udim})$ be the configuration vector of length $\udim$ such that the $i$-th element of $U_n$ is the number of balls of colour $i$ at time $n$.  Let $U_0$ be the initial vector of colour configuration, then at  time $n\geq 1$, a ball is selected uniformly at random from the urn and if the colour of the selected ball is $i$ then the ball is replaced along with $R_{i,j}$ many balls of colour $j$,  for every $1 \leq j \leq \udim$. The dynamics of the urn configuration depends on its initial configuration $U_0$ and the $\udim \times \udim$ replacement matrix $R = (R_{i,j})_{1\leq i,j\leq \udim}$.

We  study the limiting properties of the numbers of cherries and pitchforks via an equivalent urn process. Towards this, we use six different colours and assign one colour to each type of edges of a tree $T$ in the following scheme introduced in~\cite{Paper1}:  colour $1$ for all pendant edges of a cherry in a pitchfork; colour $2$ for pendant edges of an essential cherry; colour $3$ for pendant edges in a pitchfork but not in any cherry; colour $4$ for pendant edges in neither a cherry nor a pitchfork; colour $5$ for the internal edge of an essential cherry (i.e., those adjacent to colour $2$ edges), and colour $6$ for all other (necessarily internal) edges. 
See Fig.~\ref{fig:example} for an illustration of the scheme. For $1\le i \le 6$, let $E_i(T)$ be the set of edges of color $i$ in $T$.

Consider an urn with colour configuration at time $n$ as
$U_n = (U_{n,1},\dots, U_{n,6})$, where $U_{n,i}$ denotes the number of edges of colour $i$ in the tree at time $n$, which has precisely $n+2$ leaves. Then $U_0 = (0,2,0,0,1,0)$, since at the initial time step $(n=0)$ the tree $T_2$ is an essential cherry which has two pendant edges and one interior edge; see $T_2$ in Fig.~\ref{fig:example}. Based on the colouring scheme of the edges, at any time $n\geq 0$, we have
\begin{equation}\label{urn-pf-ch}
(A_{n+2}, C_{n+2}) = \frac{1}{2} \left(U_{n,1}, U_{n,1}+U_{n,2}\right),
\end{equation}
where $A_{n+2}$ and $C_{n+2}$ are the numbers of pitchforks and cherries in $T_{n+2}$, respectively.
Under the alpha tree model, the dynamics of the corresponding urn process evolves according to the following replacement matrix
\[ R = \begin{bmatrix}
0&0&0&1&0&1\\
2&-2&1&0&-1&2\\
-2&4&-1&0&2&-1\\
0&2&0&-1&1&0\\
2&-2&1&0&-1&2\\
0&0&0&1&0&1
\end{bmatrix}. \]
Let $e_i$, $1\le i \le 6$, denote a $6$-vector in which the $i$-th component is $1$ and $0$ elsewhere; and   $\chi_n$   the random vector taking value $e_i$ if, at time $n$, speciation happens at an edge with type  $i$.
Thus,  we have the following recursion
\[U_n= U_{n-1} +\chi_n R, \qquad n\geq 1, \]
where
\begin{equation}\label{SelecProb}
P(\chi_n =e_i|\FF_{n-1}) \propto \begin{cases}
(1-\al) U_{n-1,i}, &\text{ for } i \in\{1,2,3,4\}, \\[1ex]
\al \,U_{n-1,i}, & \text{ for } i\in \{5,6\}.
\end{cases}
\end{equation}
Observe that the process $(U_n)_{n \ge 0}$, which  describes the dynamics of the numbers of cherries and pitchfork,
is a {\em nonuniform urn model} since the balls are not selected uniformly at random from the urn, which is different from the classical {\em uniform} urn models in which  the balls are selected uniformly at random from the urn (see, e.g.~\cite[Chapter 7]{hofri2019algorithmics}).

We end this subsection with the following observation relating the edge color scheme with 
the number of  pitchforks and that of cherries, which follows directly from the replacement matrix $R$ (see also~\cite[Section 2]{WuChoi16}).

\begin{lemma}	\label{lem:edge-set}
	Suppose that $T$ is a phylogenetic tree with $n\ge 2$ leaves. 
	Then
	we have
	$|E_{3}(T)|=\pf(T)$,
	$|E_{2}(T)\cup E_5(T)|=3(\ch(T)-\pf(T))$,
	$|E_{4}(T)|=n-\pf(T)-2\ch(T)$,
	and
	$|E_{1}(T)\cup E_6(T)|=n-1+3\pf(T)-\ch(T)$.
Furthermore, suppose that $e$ is an edge in $T$ and $T'=T[e]$. Then
		we have
		{
			\begin{equation*}
			\pf(T') = \begin{cases}
			\pf(T)-1 & \text{if } e\in E_3(T),\\
			\pf(T)+1 & \text{if } e\in E_2(T)\cup E_5(T),\\
			\pf(T) & \text{otherwise}; \\
			\end{cases}  \quad
			\mbox{and} ~~
			\ch(T') = \begin{cases}
			\ch(T)+1 & \text{if } e \in E_3(T)\cup E_4(T),\\
						{} & \\
			\ch(T) & \text{otherwise}. \\

			\end{cases}
			\end{equation*}
		}

\end{lemma}

\subsection{Limiting theorems on uniform urn models}
\label{limit:uniform:urn}
In this subsection, we recall the strong law of large numbers and the central limit theorem on a version of uniform urn models developed in~\cite{Paper1}, which will be related to the nonuniform urn process in Subsection~\ref{subsect:urn:tree} later using the urn coupling idea  in \cite{Kaur2018}.

For the classical uniform urn models, it has been shown (see \cite{BaiHu2005}) that the random process $U_n/n$ converges almost surely to the left eigenvector of $R$ corresponding to the maximal eigenvalue and  asymptotic normality holds with a known limiting variance matrix under certain assumptions on $R$.  Standard assumptions made in the urn model theory are  that the replacement matrix is  irreducible with a constant row sum and all the off-diagonal elements are non-negative~(see, e.g.~\cite{Hosam2009}).  In  \cite{Paper1}, we extend this to the case when off-diagonal elements of a replacement matrix can be negative satisfying  the following set of assumptions {\bf (A1)--(A4)}, which was slightly rephrased from~\cite{Paper1}. Let $\text{diag}(a_1, \dots , a_\udim)$ denote the diagonal matrix whose diagonal elements are $a_1,\dots , a_\udim$. 

\noindent
{\bf (A1):} {\em Tenable:} It is always possible to draw balls and follow the replacement rule.
\\
{\bf (A2):} {\em Small:} 
All eigenvalues of $R$ are real; the maximal eigenvalue $\lambda_1$, called the {\em principal eigenvalue} is positive with $\lambda_1>2\lambda$ holds for all other eigenvalues $\lambda$ of $R$. \\
{\bf (A3):} {\em Strictly balanced:}  	The column vector $\mathbf{u}_1=(1,1,\dots,1)^\top$,  is a right eigenvector 
of $R$ corresponding to $\lambda_1$;   and it has a principal left eigenvector $\bf{v}_1$ (i.e., the 
left eigenvectors corresponding to $\lambda_1$) that is also a  probability  vector.
\\
{\bf (A4):} {\em Diagonalisable:} There exists an invertible matrix $V$ with real entries whose first row is $\bf{v}_1$ such that the first column of $V^{-1}$ is $\mathbf{u}_1$ and 
\begin{equation}
\label{eq:R:diagonal}
VRV^{-1} = \text{diag}(\lambda_1, \lambda_2,\dots, \lambda_\udim) =: \Lambda,
\end{equation}
where $\lambda_1> \lambda_2 \ge \dots \ge \lambda_\udim$ are  eigenvalues of $R$.

\medskip

Let $\mathcal{N}(\mathbf{0}, \Sigma)$ be the multivariate normal distribution with mean vector $\mathbf{0}=(0,\dots,0)$ and covariance matrix $\Sigma$. Then we have the following result from~\cite[Theorems 1 \& 2]{Paper1}, which can also be alternatively derived from~\cite[Theorems 3.21 \& 3.22 and Remark 4.2]{Janson2004}. 

\begin{theorem}
	\label{thm1&2:paper1}
	Under assumptions {\em \bf{(A1)--(A4)}}, we have
	\begin{equation}
	\label{eq:asconv:urn}
	(n\pev)^{-1} U_n \cas \mathbf{v}_1
	~\quad~\mbox{and}~\quad~ 
	n^{-1/2} (U_n - n\pev \mathbf{v}_1) \wconv \NN(\mathbf{0}, \Sigma),
	\end{equation}
	where $\pev$ is the principal eigenvalue and $\mathbf{v}_1$ is the principal left eigenvector of $R$, and 
	\begin{equation}
	\Sigma = \sum_{i,j=2}^\udim \frac{\pev \lambda_i \lambda _j \bu_i^\top \mbox{\em diag}(\mathbf{v}_1) \bu_j }{\pev-\lambda_i -\lambda_j} \mathbf{v}_i^\top \mathbf{v}_j,
	\end{equation}
	where  $\bv_j$ is the $j$-th row of $V$  and $\bu_j$ the $j$-th column of $V^{-1}$ for $2\le j \le \udim$.
\end{theorem}

\section{Limit Theorems for the Joint Distribution}
\label{sec:limiting}

In this section, we present the strong law of large numbers and the central limit theorems on the joint distribution of the number of cherries and the number  of pitchforks under the Ford model.

\subsection{Main convergence results}

For later use, we consider the following polynomials in $\alpha$: 
\begin{equation}
\label{eq:main:cov:poly}
\begin{matrix*}[l]
\gp_1 = 8\al^3-32\al^2+45\al-23, & \quad\quad\quad \gp_4 = 8\al^3-40\al^2+37\al+13, \\
\gp_2 = 40\al^3-164\al^2+221\al-97, &  \quad\quad\quad  \gp_5 = 40\al^3-112\al^2-31\al+181, \\
\gp_3 = 56\al^3-248\al^2+367\al-181, &  \quad\quad\quad  \gp_6 = 8\al^3+4\al^2-71\al+71; 
\end{matrix*}
\end{equation}
and for simplicity of notation, we do not indicate the $\phi_i$'s as  functions of $\alpha$. Moreover, it can be verified directly that $\phi_1, \phi_2, \phi_3 <0$ and $\phi_4, \phi_5, \phi_6 > 0$ for $\alpha \in (0, 1).$
Then, we have the following result on the joint asymptotic properties of the urn model process associated with the $\alpha$-tree model.

\begin{theorem}
	\label{thm:urn:edge}
	Suppose $(U_n)_{n\geq 0}$  is the urn process associated with the Ford model with parameter $\alpha\in(0,1)$. 
	Then, 
	\begin{equation}
	\label{eq:conv:urn}
	\frac{U_n}{n} \cas 
	\mathbf{v}	
	~\quad~\mbox{and}~\quad~
	\frac{ U_n- n\mathbf{v} }{\sqrt{n}} \wconv \mathcal{N}\left (\bzero,\Sigma \right),
	\end{equation}
	as $n \to \infty$, where 
	\begin{equation}\label{main:Leftev1}
	\bv= \frac{1}{2(3-2\al)} \left(2(1-\al), \,2(1-\al),\, (1-\al),\, 1+\al,\, 1-\al, \,5-3\al \right)
	\end{equation}
and with the polynomials $\gp_1,\dots,\gp_6$ defined in~\eqref{eq:main:cov:poly},
\begin{equation}
\label{eq:sigma:def}
\Sigma=\frac{1-\al}{4(3-2\al)^2(5-4\al)(7-4\al)} 
\begin{bmatrix*}[r]
-12\gp_1 & 4\gp_2 & -6\gp_1  & -2\gp_4 & 2 \gp_2 & -2 \gp_2 \\
4\gp_2 & -4\gp_3 & 2\gp_2  & -2\gp_6 & -2 \gp_3 & 2 \gp_3   \\
-6\gp_1 & 2\gp_2  & -3\gp_1  & -\gp_4 &  \gp_2 & - \gp_2  \\
-2\gp_4  &  -2\gp_6  &  -\gp_4 & \gp_5 & - \gp_6 &  \gp_6  \\
2 \gp_2   &  -2 \gp_3  &  \gp_2& - \gp_6 & - \gp_3 &  \gp_3  \\
-2 \gp_2    &  2 \gp_3   &  - \gp_2 & \gp_6  &   \gp_3 & -\gp_3  
\end{bmatrix*}.
\end{equation}
\end{theorem}
The proof of Theorem~\ref{thm:urn:edge} is given at  the end of this section.

\begin{remark} 
\label{rem3.1} For later use, here we present the limiting results on the urn model using a scaling factor relating to the time $n$ (which is motivated by noting that the number of leaves in the tree  at time $n$ is $n+2$). However, the results can be readily rephrased using the proportion of color balls in the urn process.

\end{remark}

\begin{remark} 
	Using the approach outlined in~\cite{Paper1},  Theorem \ref{thm:urn:edge} 	continues to hold for  the unrooted $\alpha$-tree models. 
\end{remark}

With Theorem~\ref{thm:urn:edge}, we are ready to present one of our main results in this paper 
concerning limit theorems on the joint distribution of the number of cherries $C_n$ and the number of pitchforks $A_n$
under the Ford model.

\begin{theorem} \label{Thm:Convg-ChPh}
Under the Ford model with parameter $\al \in [0,1]$, we have
	\[\frac{1}{n} (A_n,C_n) \cas (\nu, \mu) := \frac{1-\al}{2(3-2\al)} (1,2),\]
	and
	\[ \frac{(A_n, C_n) -n (\nu, \mu) }{\sqrt{n}} \wconv \mathcal{N}\big((0,0), S \big), \]
	where \begin{equation}
	S= \begin{bmatrix} \tau^2 & \rho \\ \rho & \si^2 \\ \end{bmatrix} =
	\frac{1-\al}{(3- 2\al)^2 (5-4\al) } \begin{bmatrix}
	\frac{-24 \al^3 +96\al^2 -135\al +69 }{4(7-4\al)} & \frac{-(2-\al) (1-2\al)}{2} \\[1ex]
	\frac{-(2-\al) (1-2\al)}{2} & 2-\al
	\end{bmatrix}.
	\end{equation}
\end{theorem}
\medskip
 \begin{remark} 
We consider special cases of $\alpha$-tree model, which are commonly studied in phylogenetics. The first two have been established in ~\cite{Paper1}.

\begin{enumerate}
	\item The uniform model corresponds to  $\al=1/2$, where all edges, internal or leaf, are selected with equal weight and the limit results hold with
	\[(\nu, \mu) = \frac{1}{8}(1, 2)
	\quad \text{and } \quad   \begin{bmatrix} \tau^2 & \rho \\ \rho & \si^2 \\ \end{bmatrix}   =
	\frac{1}{64}\begin{bmatrix} 3&0\\0&4\end{bmatrix}. \]
	\item The Yule model corresponds to  $\al =0$, where only leaf edges are selected with equal weight and the limit results hold with
	\[ (\nu, \mu) = \frac{1}{6} (1,2)
	\quad \text{and } \quad  \begin{bmatrix}
	\tau^2 & \rho \\
	\rho & \si^2 \\
	\end{bmatrix}  =  \frac{1}{45}  \begin{bmatrix} 69/28  &-1  \\  -1 & 2 \end{bmatrix}.  \]
	
\item The  Comb model corresponds to $\al =1$,  a degenerate case. It is easy to see that  $(\nu, \mu) = (0,0) $ and $\tau^2= \rho = \si^2  = 0$.
	
\end{enumerate}
 \end{remark}

\begin{proof}[Proof of Theorem~\ref{Thm:Convg-ChPh}]
	First note that the case $\al=1$ reduces to a degenerate case of Comb model and therefore we only consider $\al\in[0,1)$. The limiting results for the case $\al=0$ has been obtained in \cite{Paper1}, which agree with the above results when $\al=0$. Thus, it is enough to prove the result for $\al\in (0,1)$.
	
	 By~\eqref{urn-pf-ch}, we have 	 
	$ (A_n, C_n) = U_n \conv $
	with
	\begin{equation}\label{Def-T}
	\conv^\top = \frac{1}{2}\begin{bmatrix} 1 &0&0&0&0&0\\
	1&1&0&0&0&0
	\end{bmatrix}.
	\end{equation}		
Since
	\begin{equation}
	\frac{U_n }{n} \cas  \bv  = \frac{1}{2(3-2\al)} \big(2(1-\al), 2(1-\al), 1-\al, 1+\al, 1-\al, 5-3\al \big),
	\end{equation}
	using the relation from equation \eqref{urn-pf-ch} we get
	\[\frac{1}{n} (A_n,C_n) =\left(\frac{U_n}{n} \right) \conv \cas \bv\,\conv = \frac{1-\al}{2(3-2\al)} (1,2). \]
	This concludes the proof of the almost sure convergence. We now prove the central limit theorem and obtain the expression for the limiting variance matrix.

Denoting the $(i,j)$-entry in the covariance matrix $\Sigma$ of~\eqref{eq:sigma:def} by $\sigma_{i,j}$ for $1\le i,j\le 6$, we   consider the matrix 
	\begin{align*}
	S&= \conv^\top  \Sigma \conv
	=\frac{1}{4} \begin{bmatrix}
	\sigma_{1,1} & \sigma_{1,1}+\sigma_{1,2}\\
	\sigma_{1,1}+\sigma_{2,1}& \sigma_{1,1}+\sigma_{2,1}+\sigma_{1,2}+\sigma_{2,2}
	\end{bmatrix} \\[1ex]
	&=\frac{1-\al}{16(3-2\al)^2(5-4\al)(7-4\al)} 
	\begin{bmatrix}
	-12\gp_1 & -12\gp_1+4\gp_2\\
	-12\gp_1+4\gp_2 & -12\gp_1+8\gp_2-4\gp_3
	\end{bmatrix} \\[1ex]
		&=\frac{1-\al}{(3- 2\al)^2 (5-4\al) } \begin{bmatrix}
	\frac{-24 \al^3 +96\al^2 -135\al +69 }{4(7-4\al)} & \frac{-(2-\al) (1-2\al)}{2} \\[2ex]
	\frac{-(2-\al) (1-2\al)}{2} & 2-\al
	\end{bmatrix}.
	\end{align*}
Since  $(A_n, C_n) = U_n \conv$, where $\conv$ is as defined in \eqref{Def-T}, we get
	\[ \frac{(A_n, C_n) -n (\nu, \mu) }{\sqrt{n}} = \frac{1}{\sqrt{n}} \left(U_n -n \bv \right)\conv \wconv  \mathcal{N}\left(\mathbf{0}, \conv^\top  \Sigma \conv \right)
	=\mathcal{N}\left(\mathbf{0}, S \right)
	.\]
 This completes the proof. 
\end{proof}

\begin{figure}[H] 
	\label{fig:limit:cov}
	\centering
	\includegraphics[width=0.8\linewidth,height=0.38\textheight]{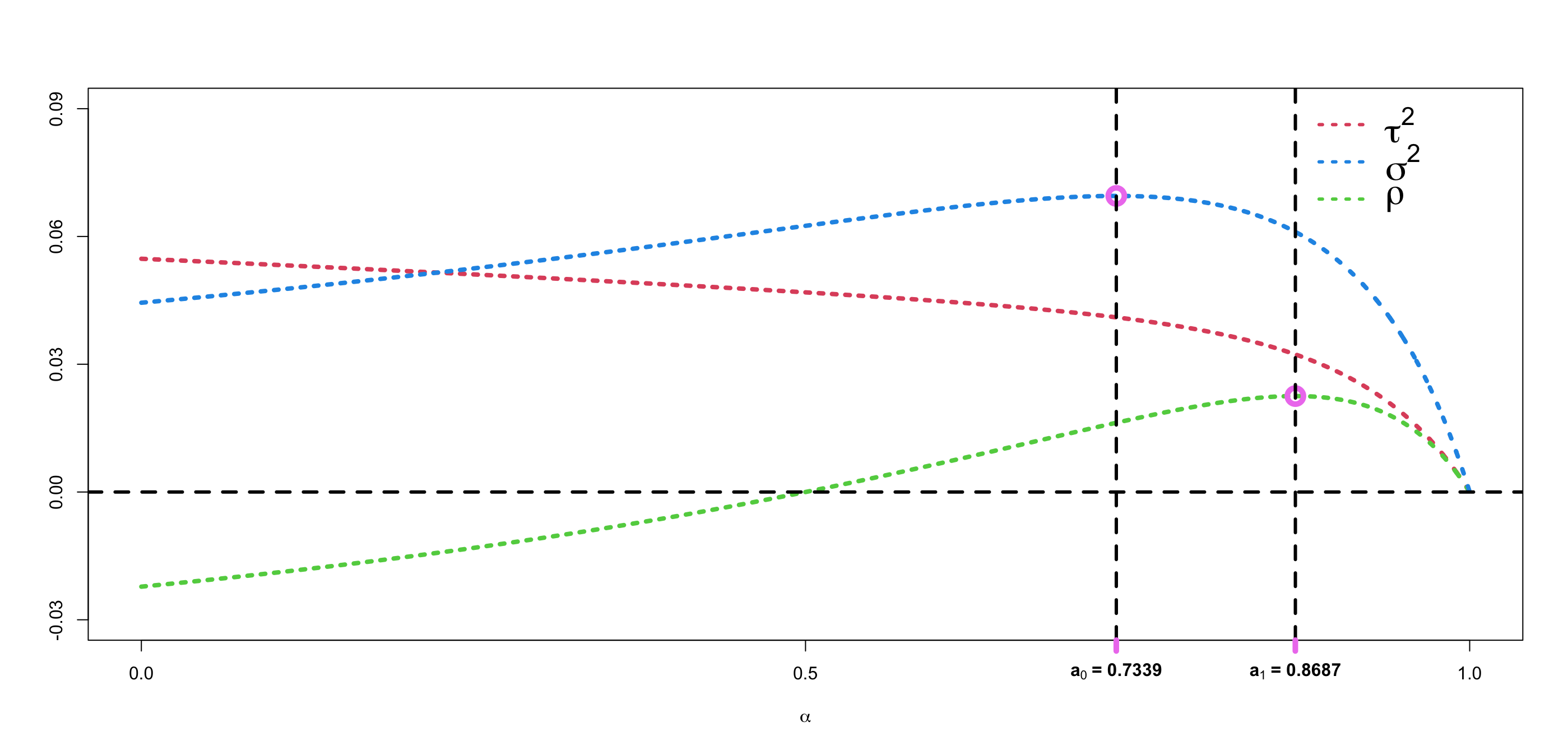}
	\caption{Plot of the limiting  variances and covariance of the joint distribution of cherries and pitchforks with respect to the parameter $\al$ under the Ford model.}
\end{figure}

We end this subsection with the following results on the behaviour of the first and second moments of the limiting joint distribution of cherries and pitchforks in the parameter region, as indicated by their plots in Figure~\ref{fig:limit:cov}.

\begin{corollary}
	\begin{enumerate}
	\renewcommand{\theenumi}{\roman{enumi}}
		\item For $0 \le \al <1$,  $A_n/C_n  \cas 1/2$ as $n\to \infty$. That is, the number of pitchforks is asymptotically equal to the number of essential cherries.
		
		\item $ A_n/n \cas  \frac{1-\al}{2(3-2\al)}, $ this limit  decreases strictly from $1/6$ to $0$, as $\al$ increases from $0$ to $1$.

		\item  The limiting variance of $A_n/\sqrt{n}$,  $\tau^2$, decreases strictly from $23/420$ to $0$, as $\al$ increases from $0$ to $1$.
		\item The limiting variance of $C_n/\sqrt{n} $, $\si^2$, increases strictly from $2/45$ 
		to $0.0695$  over $(0, a_0)$ and decreases from $0.0695$ to $0$ over $(a_0, 1)$, where $a_0 =0.7339$, 
		 the unique root of $19-48\al+36\al^2-8\al^3 =0$ in $(0,1)$.
		
		\item  The limiting covariance of $A_n/\sqrt{n}$ and $C_n/\sqrt{n}$ changes sign from negative to positive at $\al =1/2$. Specifically, it increases from $-1/45$  
		to $0.0225$  
		over $(0, a_1)$ and decreases from $0.0225$ 
		over $(a_1, 1)$, where $a_1=0.8688$, 
		the unique root of $-24\al^4+160\al^3 -370\al^2 +358\al -123=0$ in $(0,1)$.
	\end{enumerate}
\end{corollary}

\subsection{A uniform urn model derived from $U_n$}

For $\al\in (0,1)$, consider the diagonal $6\times 6$ matrix
$
\proj_\al=\text{diag}(1-\al,1-\al,1-\al,1-\al,\al,\al)
$ 
and \[ \tU_n := U_n\proj_\al = \left((1-\al) U_{n,1}, \dots, (1-\al) U_{n,4}, \al U_{n,5}, \al U_{n,6}\right). \]
Clearly, there is a one to one correspondence between $U_n$ and $\tU_n= U_n\proj_\al$ for $\al\in (0,1)$ and therefore it is sufficient to obtain the limiting results for the urn process $\tU_n$. Note that the off-diagonal elements of the replacement matrix $R_\al$ are not all  non-negative, therefore we will use the limit results from \cite{Paper1} to  obtain the convergence results for  the urn process $\tU_n$.

\begin{theorem}\label{Thm:Urn}
Suppose $\alpha\in(0,1)$. Then 	$(\tU_n)_{n\geq 0}$ is an uniform urn process with replacement matrix $ R_\al = R\proj_\al$ and
\begin{equation}
\label{eq:cas:nonuniform}
\frac{\tU_n }{n} \cas \widetilde{\bv}_1,
\end{equation}
where \begin{equation}\label{Leftev1}
\widetilde{\bv}_1= \frac{1}{2(3-2\al)} \big(2(1-\al)^2, 2(1-\al)^2, (1-\al)^2, 1-\al^2, \al(1-\al), \al(5-3\al) \big)
\end{equation}
is the normalized left eigenvector of $R_\al$ corresponding to the largest eigenvalue $\lambda_1=1$. Furthermore,
\begin{equation}
\label{eq:cweak:nonuniform}
\frac{\tU_n -n \widetilde{\bv}_1}{\sqrt{n}} \wconv  \mathcal{N}(\mathbf{0}, \widetilde{\Sigma}),
\end{equation}
with the polynomials $\gp_1,\dots,\gp_6$ defined in~\eqref{eq:main:cov:poly} and $\nal=1-\al$,
\begin{equation}
\label{eq:covariance:nonuniform}
\widetilde{\Sigma}=
\frac{\nal}{4(3-2\al)^2(5-4\al)(7-4\al)} 
\begin{bmatrix*}[r]
-12\nal^2\gp_1 & 4\nal^2\gp_2 & -6\nal^2\gp_1  & -2\nal^2\gp_4 & 2\al\nal \gp_2 & -2\al\nal \gp_2 \\
4\nal^2\gp_2 & -4\nal^2\gp_3 & 2\nal^2\gp_2  & -2\nal^2\gp_6 & -2\al\nal \gp_3 & 2\al\nal \gp_3   \\
-6\nal^2\gp_1 & 2\nal^2\gp_2  & -3\nal^2\gp_1  & -\nal^2\gp_4 & \al\nal \gp_2 & -\al\nal \gp_2  \\
-2\nal^2\gp_4  &  -2\nal^2\gp_6  &  -\nal^2\gp_4 & \nal^2\gp_5 & -\al\nal \gp_6 & \al\nal \gp_6  \\
2\al\nal \gp_2   &  -2\al\nal \gp_3  &  \al\nal \gp_2& -\al\nal \gp_6 & -\al^2 \gp_3 & \al^2 \gp_3  \\
-2\al\nal \gp_2    &  2\al\nal \gp_3   &  -\al\nal \gp_2 & \al\nal \gp_6  &  \al^2 \gp_3 & -\al^2 \gp_3  
\end{bmatrix*}.
\end{equation}

\end{theorem}
\begin{proof}[Proof of Theorem \ref{Thm:Urn}]
	First, observe that at any time $n$, there are $n+2$ pendant edges and $n+1$ internal edges in a rooted tree. That is, 
\[ U_{n,1}+ U_{n,2}+ U_{n,3}+ U_{n,4}= n+2 \quad \text{and } \quad U_{n,5}+ U_{n,6} = n+1.\]
This gives
\begin{align*}
\|\tU_n\|_1
= (1-\al) \sum_{j=1}^4 U_{n,j}+ \al  \sum_{j=5}^6 U_{n,j}  
= (1-\al) (n+2) + \al (n+1) = n+2 -\al.
\end{align*}
Therefore, from \eqref{SelecProb} we get,
\[ \Ebold[\chi_{n}| \FF_{n-1}] = \dfrac{U_{n-1} \proj_{\al}}{\|U_{n-1}\proj_{\al}\|_1}= \dfrac{U_{n-1} \proj_{\al}}{ n+1 -\al},
\]
and
\begin{align*}
\Ebold[U_{n}|\FF_{n-1}]
= U_{n-1} + \Ebold[\chi_{n}|\FF_{n-1}] R
= U_{n-1} + \dfrac{1}{n+1 -\al} U_{n-1} \proj_\al R.
\end{align*}
Multiplying both sides by $\proj_\al$,  we get

\[\Ebold[\tU_{n}|\FF_{n-1}] = \tU_{n-1} + \left(\dfrac{1}{\|\tU_{n-1} \|_1} \tU_{n-1} \right) R \proj_\al. \]
Hence, $(\tU_n)_{n\geq 0}$ is a classical uniform urn model with replacement matrix $ R_\al = R\proj_\al$.

Note that  {\bf (A1)} holds because the general Ford's dynamics on a rooted tree is well defined at every time $n$, thus the corresponding urn model satisfies the assumption of tenability. That is, it is always possible to draw balls  without getting stuck  with the replacement rule.
Note that $R_\al$ is diagonalisable as 
\[V R_\al V^{-1}=\Lambda \]
holds with  $ \Lambda = \text{diag} \big(1,0,0,0, -2(1-\al),-(3-2\al)\big)$,
\begin{equation} \label{Reigen:AalphaR}
V^{-1}= \begin{bmatrix}
1& \frac{1}{\nal}&0&0&1& 1-\al\\[1ex]
1& 0&\frac{1}{\nal}&0&1& 3-\al \\[1ex]
1& \frac{-2}{\nal}&0& \frac{3}{\nal}&\frac{-(2-\al)}{\nal}&-5+\al \\[1ex]
1 &0&0&\frac{1}{\nal}&\frac{-(2-\al)}{\nal}&-3+\al \\[1ex]
1 & 0&\frac{-2}{\al}&\frac{1}{\al}&1&3-\al \\[1ex]
1 & 0&0&\frac{-1}{\al} & 1& 1-\al
\end{bmatrix}
\end{equation}

and 
	\setlength{\arraycolsep}{2.5pt}
\medmuskip = 1mu
\begin{align}
V =\frac{1}{2(3-2\al) }
\begin{bmatrix}\label{Leigen:AalphaR:1}
2\nal^2 & 2\nal^2 & \nal^2 & (1+\al)\nal & \al\nal& \al(5-3\al)\\[1ex]
2\nal(1+\al-\al^2) & 2\nal^3 &-(2-\al)\nal^2 & (2-\al)\nal^2 & -\al\nal^2& -\al\nal(5-3\al) \\[1ex]
2\al\nal^2 & 2\al(2-\al)\nal & \al\nal^2 &-\al\nal^2 & -\al(3- \al)\nal & -3\al\nal^2 \\[1ex]
2\al(2-\al)\nal & 2\al\nal^2 & \al(2-\al)\nal &-\al(2-\al)\nal & \al^2 \nal & -3\al(2-\al)\nal \\[1ex]
2(2-\al)\nal & -2\nal^2 & (2-\al)\nal & -(4-\al)\nal & -\al\nal &\al\nal\\[1ex]
-2\nal & 2\nal & -\nal & \nal & \al & -\al
\end{bmatrix}.
\end{align}
Therefore, $R$ satisfies condition {\bf (A4)}. Next, {\bf (A2)} holds because $R_\al$ has eigenvalues 
$$1,\quad 0,\quad 0, \quad 0, \quad -2(1-\al),\quad -(3-2\al)$$ which are all real. The maximal eigenvalue $\lambda_1=1$ is positive with $\lambda_1>2\lambda$ holds for all other eigenvalues $\lambda$ of $R_\al$.
Furthermore, put $\mathbf{u}_i={V^{-1}}\mathbf{e}^\top_i$ and $\mathbf{v}_i=\mathbf{e}_iV$ for $1\le i \le 4$. Then {\bf (A3)} follows by noting that $\mathbf{u}_1=(1,1,1,1,1,1)^\top$ is the principal right eigenvector, and
\[\widetilde{\bv}_1= \frac{1}{2(3-2\al)} \big(2(1-\al)^2, 2(1-\al)^2, (1-\al)^2, 1-\al^2, \al(1-\al), \al(5-3\al) \big)\]  is the principal left eigenvector.

Since all the assumptions {\bf (A1)--(A4)} are satisfied  by the replacement matrix $R_\al$, by Theorem~\ref{thm1&2:paper1}, 
\eqref{eq:cas:nonuniform} holds. Furthermore, since
\begin{equation} 
\widetilde{\Sigma} = \sum_{i,j=2}^6 \frac{ \lambda_i \lambda _j \bu_i^\top \mbox{diag}(\mathbf{v}_1) \bu_j }{1-\lambda_i -\lambda_j} \mathbf{v}_i^\top \mathbf{v}_j,
\end{equation}
by~\eqref{eq:cas:nonuniform} it follows that~\eqref{eq:cweak:nonuniform} holds.
\end{proof}

\subsection{Proof of Theorem \ref{thm:urn:edge}}

\begin{proof}
	\noindent
	Observe that  $\sum_{i=1}^6U_{n,i} = 3+2n$ (since $2$ balls are added into the urn at every time point), thus the vector of color proportions is $U_n /(3+2n)$.
	Since $\al \in (0,1)$, it follows that  $\proj_\al$ is invertible and its inverse is 
	$$
	\proj_\al^{-1}=\frac{1}{\al(1-\al)}\mbox{diag}(\al,\al,\al,\al,1-\al,1-\al),
	$$
	which is also a diagonal matrix, and so $(\proj_\al^{-1})^\top=\proj_\al^{-1}$.
	Note that  we have $U_n = \tU_n \proj_\al^{-1}$ and consider 
	$$
	{\bv} = \widetilde{\bv}_1 (\proj_\al)^{-1} = \frac{1}{2(3-2\al)} \big(2(1-\al), 2(1-\al), 1-\al, 1+\al, 1-\al, 5-3\al \big).
	$$
	Since $ \dfrac{\tU_n }{n} \cas \widetilde{\bv}_1$ holds in view of~\eqref{eq:cas:nonuniform} in Theorem \ref{Thm:Urn}, 
	\begin{equation}
	\frac{U_n }{n} \cas {\bv},
	\end{equation}
	which concludes the proof of the almost sure convergence in~\eqref{eq:conv:urn}. 
	
	Consider the covariance matrix $\widetilde{\Sigma}$ for $ \tU_n$ as stated in~\eqref{eq:covariance:nonuniform}, then by straightforward calculation we have 
	$$
	\Sigma=(\proj_\al^{-1})^\top \widetilde{\Sigma} \proj_\al^{-1}=\proj_\al^{-1}  \widetilde{\Sigma} \proj_\al^{-1}. 
	$$
	Therefore, since
	\[\frac{\tU_n -n \widetilde{\bv}_1}{\sqrt{n}} \wconv \mathcal{N} (\mathbf{0}, \widetilde{\Sigma} ) \]
	in view of Theorem \ref{Thm:Urn}, we get
	\[\frac{U_n -n \bv }{\sqrt{n}} \wconv \mathcal{N} \big(\mathbf{0}, (\proj_\al^{-1})^\top \widetilde{\Sigma}\, \proj_\al^{-1}\big)=\mathcal{N} (\mathbf{0},  \Sigma). \]
	This completes the proof.
\end{proof}


\section{Exact Distributions }
\label{sec:exact}

In this section, we present recursion formulas for exact computation of  the joint distributions of cherries and pitchforks, their means, variances and covariance for fixed $n$ under the Ford model.

We start with the following result on the exact computation of the  joint probability mass function (pmf)  of $A_n$ and $C_n$, 
which can be regarded as a generalization of the previous results on the Yule model (e.g. when $\alpha=0$~\cite[Theorem 1]{WuChoi16}) and  the uniform model (e.g. $\alpha=1/2$~\cite[Theorem 4]{WuChoi16}).
A related result for unrooted  trees is presented in~\cite{CTW19}.

\begin{theorem} \label{jointpmf}
For $n \ge 3$, $0 \le  a\le n/3$ and $1\le b\le n/2$, under the Ford model with parameter $\alpha\in [0,1]$  we have
\begin{eqnarray*}
&& \pp(A_{n+1}=a, C_{n+1}=b) \\&=&  \frac{2a+ \al(n-a-b-1)}{n-\al} \pp(A_n=a, C_n=b)
+ \frac{(1-\al)(a+1)}{n-\al} \pp(A_n=a+1, C_n=b-1) \\
&& 
\quad
+ \frac{(2-\al)(b-a+1)}{n-\al} \pp(A_n=a-1, C_n=b) + \frac{(1-\al)(n-a-2b+2)}{n-\al} \pp(A_n=a, C_n=b-1).
\end{eqnarray*}
\end{theorem}
\begin{proof}[Proof of Theorem \ref {jointpmf}]
Fix $n> 3$, and let $T_2,\dots,T_n,T_{n+1}$ be a sequence of random trees generated by the Ford process, that is, $T_2$ contains two leaves and $T_{i+1}=T_i[e_i]$ for a random edge $e_i$ in $T_i$  chosen according to the Ford model for $2\leq i \leq n$. 
	Then we have
	\begin{align}
	\label{eq:total:yule}
	\pp(A_{n+1}=a, C_{n+1}=b) &=\pp(\pf(T_{n+1})=a, \ch(T_{n+1})=b) \notag \\	
	&\hspace{-2cm}
	=\sum_{p,q} \pp(\pf(T_{n+1})=a, \ch(T_{n+1})=b\,|\,\pf(T_n)=p,\ch(T_n)=q) \pp(\pf(T_n)=p,\ch(T_n)=q) \notag \\
	\hspace{-1cm}
	&\hspace{-2cm}
	=\sum_{p,q} \pp (\pf(T_{n+1})=a, \ch(T_{n+1})=b\,|\,\pf(T_n)=p,\ch(T_n)=q) \pp(A_{n}=p, C_{n}=q),
	\end{align}
	where the first and second equalities follow from the law of total probability, and the definition of random variables $A_n$ and $C_n$.
	
Let $e_n$ be the edge in $T_n$ chosen in the above Ford process for generating $T_{n+1}$, that is, $T_{n+1}=T_n[e_n]$. Since Lemma~\ref{lem:edge-set} implies that
	\begin{equation}
	\label{eq:yule:5}
	\pp(\pf(T_{n+1})=a, \ch(T_{n+1})=b~|~\pf(T_n)=p,\ch(T_n)=q)=0 
	\end{equation}
	for $(p,q) \not \in\{(a,b),(a+1,b-1),(a-1,b),(a,b-1)\}$,
	it suffices to consider the following four cases in the summation in (\ref{eq:total:yule}): case (i): $p=a, q=b$; case (ii): $p=a+1, q=b-1$; case (iii): $p=a-1, q=b$; and case (iv): $p=a, q=b-1$.

	First, Lemma~\ref{lem:edge-set} implies that case (i) occurs if and only if $e_n\in E_1(T_n)\cup E_6(T_n)$, and  that $E_1(T_n)\cup E_6(T_n)$ contains precisely $2\pf(T_n)$ pendent edges and $(n-1)+\pf(T_n)-\ch(T_n)$ interior edges. Therefore we have
	\begin{align}
	&\pp(\pf(T_{n+1})=a, \ch(T_{n+1})=b~|~\pf(T_n)=a,\ch(T_n)=b)  \nonumber \\
	&\quad \quad =\frac{2\pf(T_n)(1-\al)+\al(n-1+\pf(T_n)-\ch(T_n))}{n-\al}=
	\frac{2a+ \al(n-a-b-1)}{n-\al}.\label{eq:yule:1}
	\end{align}

	Similarly, 
	case (ii)
	occurs if and only if $e_n\in E_3(T_n)$, which contains  $\pf(T_n)$ pendent edges and no interior edges. Therefore we have
	\begin{align}
	\label{eq:yule:2}
	\pp(\pf(T_{n+1})=a, \ch(T_{n+1})=b~|~\pf(T_n)=a+1,\ch(T_n)=b-1)
	=\frac{(a+1)(1-\al)}{n-\al}.
	\end{align}
	
	Next, 
	case (iii) occurs precisely when $e_n\in E_2(T_n)\cup E_5(T_n)$, which contains  $2(\pf(T_n)-\ch(T_n))$ pendent edges and $\pf(T_n)-\ch(T_n)$ interior edges. Thus
	\begin{align}
	&\pp(\pf(T_{n+1})=a, \ch(T_{n+1})=b~|~\pf(T_n)=a-1,\ch(T_n)=b) \nonumber \\
	&	\quad \quad
	=\frac{2(a-1-b)(1-\al)+\al(a-1-b)}{n-\al}=\frac{(2-\al)(b-a+1)}{n-\al}.\label{eq:yule:3}
	\end{align}

	Finally,  
	case (iv) occurs
	if and only if $e_n$ is  in $E_4(T_n)$,  which contains precisely $n-\pf(T_n)-2\ch(T_n)$ pendent edges and no interior edges. Hence, 
	\begin{equation}
	\label{eq:yule:4}
	\pp(\pf(T_{n+1})=a, \ch(T_{n+1}=b)~|~\pf(T_n)=a,\ch(T_n)=b-1)
=\frac{(1-\al)(n-a-2b+2)}{n-\al}. 
	\end{equation}
	
	Substituting Eq.~\eqref{eq:yule:1}--\eqref{eq:yule:4} into Eq.~\eqref{eq:total:yule} completes the proof of the theorem.
	\end{proof}

To study the moments of $A_n$ and $C_n$, we present below a functional recursion form of Theorem~\ref{jointpmf}, whose proof is straightforward and hence omitted here.

\begin{proposition}
\label{jointrr}
Let $\varphi: \mathbb{N}\times \mathbb{N} \to \mathbb{R}$ be an arbitrary function. For $n \ge 3$,
under the Ford model with parameter $\alpha\in [0,1]$ we have
\begin{eqnarray*}
(n-\al) \Ebold \varphi(A_{n+1}, C_{n+1}) &=& \Ebold\bigg[ \big\{\al(n-A_n-C_n-1) +2A_n \big\} \varphi(A_n, C_n) \\
&& +(1-\al)A_n \varphi(A_n-1, C_n+1) +(2-\al) (C_n-A_n) \varphi(A_n+1, C_n) \\
&& + (1-\al)(n-A_n-2C_n) \varphi(A_n, C_n+1) \bigg ].
\end{eqnarray*}
\end{proposition}

For a fix integer $k$, consider the indicating function $I_k(x,y)$ that equals to 1 if $y=k$, and $0$ otherwise. 
Then applying Proposition~\ref{jointrr} with $\varphi(x, y)=I_k(x)$ leads to the following result on the distribution of cherries.

\begin{corollary}
\label{cherrypmf}
For integers $n \ge 3$ and $0\le k \le n/2$, under the Ford model with parameter $\alpha\in [0,1]$ we have
$$
(n-\al) \pp(C_{n+1}=k) = [(n-1) \al +2(1-\al) k ] \pp(C_{n}=k)
+(1-\al)(n-2k +2) \pp(C_{n+1}=k-1).
$$
\end{corollary}

Similarly, applying Proposition~\ref{jointrr} with appropriate functions $\varphi$ leads to the following recurrence relation on the moments of the joint distributions; 
the proof is similar to those of~\cite[Corollary 4 \& Proposition 5]{WuChoi16} and hence omitted here.

\begin{corollary}
\label{recurrence}
For $n \ge 3$, under the Ford model with parameter $\alpha\in [0,1]$ we have
\begin{eqnarray}
(n-\al) \Ebold [C_{n+1}] - (n -2+\al) \Ebold [C_n] &=& n(1-\al), \label{cherrymean}\\
(n-\al) \Ebold[A_{n+1}] - (n-3+\al) \Ebold[A_n] &=& (2-\al) \Ebold [C_n], \label{forkmean} \\
(n-\al) \Ebold [C_{n+1}^2] - (n-4+3\al) \Ebold [C_n^2] &=& 2(n-1)(1-\al) \Ebold [C_n] + n(1-\al), \label{cherry2nd}\\
(n-\al)\Ebold[A_{n+1}C_{n+1}] - (n-5+3\al) \Ebold [A_nC_n] &=& (n-1)(1-\al)\Ebold[A_n] + (2-\al) \Ebold[ C_n^2], \label{cov}\\
(n-\al) \Ebold[A_{n+1}^2]- (n-6+3\al)\Ebold[A_n^2] &=& 2(2-\al)\Ebold [A_nC_n] + (2-\al)\Ebold [C_n] - \Ebold [A_n] \label{fork2nd}
\end{eqnarray}
with initial conditions $\Ebold[A_3]= \Ebold[C_3]=\Ebold[A_3^2]=\Ebold[C_3^2]=\Ebold[A_3C_3]=1. $
\end{corollary}


\begin{remark}
Let $\mu_n = \Ebold [C_n]$ and 	$\si_{n}^2=\var(C_n)$. 
Substituting $\Ebold[C_{n}^2]=\si_{n}^2 + \mu_{n}^2$ into (\ref{cherry2nd}) and applying (\ref{cherrymean}), we obtain below a recurrence relation of the $\si_n^2$, which was also obtained in Ford's thesis (Theorem 60, \cite{Ford2006}):
\begin{eqnarray*}
(n-\al) \si_{n+1}^2 - (n-4+3\al) \si_n^2 
&=& -\frac{4(1-\al)^2}{n-\al} \mu_n^2 + \frac{2(1-\al)[(1-2\al)n +\al]}{n-\al} \mu_n
+ \frac{\al (1-\al)n(n-1)}{n-\al}. \label{cherryvar}
\end{eqnarray*}
\end{remark}

\tw{As an application of Corollary~\ref{recurrence}, in the next theorem we obtain the formulas for the mean of $A_n$ and  that of $C_n$ under the Ford model, which extends previous results on the Yule and the uniform models (see e.g.~\cite{WuChoi16} and the references therein). Note that the mean of $C_n$ as stated in Theorem~\ref{thm:firstm} (i) was first obtained in  \cite{Ford2006} and is included here for completeness.}
 
\begin{theorem}
\label{thm:firstm} Under the Ford model with parameter \tw{$\alpha\in [0,1]$},
we have
\begin{enumerate}
\renewcommand{\theenumi}{\roman{enumi}}
\item $ \Ebold [C_n]  = \dfrac{1-\al}{3-2\al} \ n + \dfrac{\al}{2(3-2\al)} +x_n, $
where $x_3 = \frac{\al}{2(3-2\al)}$ and for $n\ge 4$,
\begin{equation}
\label{Order:xn}
x_n = \frac{\al }{2(3-2\al)} \prod_{i=3}^{n-1} \frac{i-2+\al}{i-\al}= \frac{\al\Ga(3-\al)}{2(3-2\al)\Ga(1+\al)}n^{-2(1-\al)} \left(1+o(1) \right);
\end{equation}

\item $\Ebold [A_n]  = \dfrac{1-\al}{2(3-2\al)} \ n + \dfrac{\al}{2(3-2\al)} + y_n,$
where 
$y_3=\frac{1}{2}$, 
$y_4=\frac{\al(5-3\al)}{2(3-\al)(3-2\al)}$, and  for $n\ge 5$,
\begin{equation}\label{Order:yn}
y_n = \frac{\al(2n-3+\al-n\al) }{2(3-2\al)(3-\al)}  \prod_{j=4}^{n-1}\frac{j-3+\al}{j-\al }
= \frac{\al(2-\al) }{2(3-2\al)} \frac{\Gamma(3-\al) }{\Gamma(1+\al)} n^{-2(1-\al)} \left(1+o(1) \right).
\end{equation}
\end{enumerate}
\end{theorem}


\tw{The proof of Theorem~\ref{thm:firstm} and that of Theorem~\ref{thm:secondm}, which concerns higher order expansions of the second moments, are presented in Section~\ref{sec:expansion}}.

\begin{theorem}
\label{thm:secondm}
Under the Ford model with parameter 
\tw{$\alpha\in [0,1]$}
we have
\begin{enumerate}
\renewcommand{\theenumi}{\roman{enumi}}
\item $$\var (C_n) = \frac{(1-\al)(2-\al)}{(3-2\al)^2(5-4\al)} \ n -  \frac{\al (1-\al)(2-\al)}{(3-2\al)^2(5-4\al)} +\OO(n^{-2(1-\al)}).$$
\item
$$\cov(A_n, C_n)  = \frac{-(1-\al)(2-\al)(1-2\al)}{2(3-2\al)^2(5-4\al) } \ n -\frac{ \al (1-\al)(2-\al)}{(3-2\al)^2(5-4\al)}  +\OO(n^{-2(1-\al)}).$$
\item
$$ \var(A_n) =  \frac{(1-\al)(69-135\al+96\al^2-24\al^3)}{4(3-2\al)^2(5-4\al)(7-4\al)} \ n   + \frac{3\al(1-\al)(1-2\al)(5-3\al)}{4(3-2\al)^2(5-4\al)(7-4\al)}  +\OO(n^{-2(1-\al)}).$$
\end{enumerate}
\end{theorem}


Let $\rho_\al(A_n,C_n)$ be the correlation of $A_n$ and $C_n$ under the Ford model with parameter $\al\in [0,1)$, which is not defined for  $\al=1$ because in this case $A_n$ and $C_n$ are both degenerate random variables.
It is shown in~\cite[Corollaries~3~\& 5]{WuChoi16}
 that for the Yule model  $\rho_0(A_n,C_n)=-\sqrt{14/69}$ holds for $n\ge 7$, and for the uniform model $\{\rho_{1/2}(A_n,C_n)\}_{n\ge 4}$ is an increasing sequence converging to $0$. 
 Together with Theorem~\ref{thm:secondm}(ii), this leads directly to the following result which shows that $\al=1/2$ is a critical value for $\rho_\al(A_n,C_n)$: 
 when $n$ is large, 
$A_n$ and $C_n$ are negatively correlated  for $ \al \in [0, 1/2]$, which is expected, and  positively correlated for $ \al \in (1/2, 1)$, which is less expected.

\begin{corollary}
Under the Ford model, for each $0\le \al \le 1/2$ there exists a constant $n_0(\al)$ such that 
$\rho_\al(A_n,C_n)<0$ for all $n>n_0(\al)$. Furthermore, for each $1/2<\al<1$ there exists a constant $n_0(\al)$ such that 
$\rho_\al(A_n,C_n)>0$ for all $n>n_0(\al)$.
\end{corollary}



\section{Proofs of Theorems~\ref{thm:firstm} and~\ref{thm:secondm}} 
\label{sec:expansion}
In this section, we present the proofs of the two theorems, starting with the two lemmas below.

\begin{lemma}\label{Lemma1}
Let $a, b$ and $c$ be three positive real numbers with \tw{$a > b-1$.}
Given an integer $n_0\ge 2$, suppose that 
 $\{X_n\}_{n \ge n_0}$ is a sequence of real numbers satisfying the recursion
\[ X_{n+1}  = f_n X_n +g_n, \qquad n \geq n_0,\]
where  $\{f_n\}_{n \ge n_0}$ and $\{g_n\}_{n\ge n_0}$ are two sequences with  \tw{$ \prod_{i=\ell}^{n-1} \left|f_i \right| \leq c (n/\ell)^{-a} $
}and  $|g_\ell| \leq c \ell^{-b}$ for every $ \ell \geq n_0$. 
Then there exists a positive number $C$ 
such that  $|X_n | \leq C n^{1-b}$ for all $n\ge n_0$.
\end{lemma}

\begin{proof}[Proof of Lemma \ref{Lemma1}]
Since the  solution to the given recursion is given by
\[X_n =X_{n_0} \prod_{i= n_0}^{n-1} f_i + \sum_{i=n_0}^{n-1} g_i  \prod_{j=i+1}^{n-1}f_j\]
for $n > n_0$, we have
\[ |X_n| \leq  |X_{n_0}| \prod_{i= n_0}^{n-1}\left|  f_i \right| + \sum_{i=n_0}^{n-1} |g_i|\prod_{j=i+1}^{n-1} \left|  f_j  \right|.\]
Considering $C=2 \max\{c(n_0)^a |X_{n_0}|,\frac{c^2\, 2^a}{a-b+1}\}$, then the lemma follows by noticing that
\begin{equation*}
 |X_{n_0}| \prod_{i= n_0}^{n-1}  \left|f_i \right|  \leq 
c{n_0}^a |X_{n_0}| n^{-a}
\leq C n^{-a}/2
\leq Cn^{1-b}/2
\end{equation*}
\noindent
and
\begin{align*}
 \sum_{i=n_0}^{n-1} |g_i| \prod_{j=i+1}^{n-1} \left|  f_j \right|
 &\leq c  \sum_{i=n_0}^{n-1} |g_i|   \left(\frac{i+1}{n}\right)^{a}  \leq  
 \frac{c^2\, 2^a}{n^a}
 \,  \sum_{i=n_0}^{n-1} i^{\,a-b}   
 \leq  \frac{c^2\, 2^a \, n^{a-b+1}}{n^a (a-b+1)}
\leq \frac{C}{2}\, n^{1-b}.
\end{align*}
\end{proof}

\begin{lemma} \label{Lemma2}
For $\alpha\in [0,1]$ and three finite non-negative integers $l,k,m$ such that $l\geq k$ and $m\geq 1$, there exists a positive constant $K=K(\al,m)$ 
such that
\begin{equation} \label{eq:bound1}
	\prod_{i=l }^{n-1} \left|\frac{i-k+m\al}{i-\al}  \right|   \leq
	K \left(\frac{n}{l}\right)^{-k+(m+1)\al}
	~~\mbox{for all $1\le l\le n-1$.}
	\end{equation}
Furthermore, as $n\to \infty$ we have
\begin{equation}\label{Order:prod}
\prod_{i=l }^{n-1} \frac{i-k+m\al}{i-\al}   =   \frac{\Gamma(l-\al) }{\Gamma(l-k+m\al)} n^{-k+(m+1)\al} \left(1+o(1) \right).
\end{equation}
\end{lemma}

\begin{proof}[Proof of Lemma \ref{Lemma2}]
First, \eqref{eq:bound1} follows from \cite[Lemma 2]{Paper1}. To prove \eqref{Order:prod},
 note that
\begin{align}
\prod_{i=l }^{n-1} \frac{i-k+m\al}{i-\al}
& =  \prod_{i=l }^{n-1} \frac{\Gamma(i+1-k+m\al)\Gamma(i-\al) }{\Gamma(i-k+m\al)\Gamma(i+1-\al) } 
 =   \frac{\Gamma(n-k+m\al) }{\Gamma(l-k+m\al)}   \frac{\Gamma(l-\al) }{\Gamma(n-\al)} \nonumber\\
& =   \frac{\Gamma(l-\al)}{\Gamma(l-k+m\al)}  \frac{\Gamma(n+m\al)  }{\Gamma(n-\al)}   \prod_{j=1}^k \frac{1}{n-j+m\al}.
\label{Prod:Gamma}
\end{align}

By  Stirling's approximation formula, $\Ga(x) = \sqrt{2 \pi} \ x^{x-1/2} e^{-x} \left(1+o(1)\right)$, we have
\begin{eqnarray}
\frac{\Ga(n+m\al)}{\Ga(n-\al)}
&=& \frac{\sqrt{2\pi} (n+m\al)^{n+m\al-1/2} \, e^{-(n+m\al)}}{\sqrt{2\pi} (n-\al)^{n-\al-1/2} \, e^{-(n-\al)}}\left(1+o(1) \right) \nonumber \\
&=& n^{(m+1)\al} \frac{(1+ m\al/n)^{n+m\al-1/2}}{(1-\al/n)^{n-\al-1/2}} e^{-(m+1)\al} \left(1+o(1) \right) 
=n^{(m+1)\al} \left(1+o(1) \right). \label{Eq:Str}
\end{eqnarray}
Combining \eqref{Eq:Str} and
$\prod_{j=1}^k \frac{1}{n-j+m\al} =n^{-k}\left(1+o(1) \right)$, we get \eqref{Order:prod}.
\end{proof}

Next, we present the proof of the first theorem.

\begin{proof}[Proof of Theorem~\ref{thm:firstm}]
To prove part (i), we consider
\begin{equation}\label{Exp:mu}
x_n=
\Ebold[C_n]- \frac{1-\al}{3-2\al} \ n - \frac{\al}{2(3-2\al)}, \quad n \ge 3.
\end{equation}
Since $ \Ebold[C_3]=1$, we get
$x_3  
= \frac{\al}{2(3-2\al)}.$
Furthermore, substituting (\ref{Exp:mu}) into (\ref{cherrymean}) leads to
$$(n-\al) x_{n+1}- (n-2+\al)x_n=0, \quad n \ge 3,  $$
and hence
\[x_n 
=
x_3  \prod_{i=3}^{n-1} \frac{i-2+\al}{i-\al} = x_3 \frac{ \Gamma(3-\al)}{ \Ga(1+\al)} \frac{\Ga(n-2+\al)}{\Ga(n-\al)}, \quad n \ge 4. \]
Together with Lemma \ref{Lemma2}, this establishes \eqref{Order:xn}, and hence completes the proof of part (i).\\

To prove part (ii), 
 we consider
\begin{equation}\label{Exp:nu}
y_n=\Ebold[A_n] - \frac{1-\al}{2(3-2\al)} \ n - \frac{\al}{2(3-2\al)} 
\end{equation}
for $n\ge 3$. 
Then $y_3=1/2$. Furthermore, substituting \eqref{Exp:nu} and \eqref{Exp:mu} into \eqref{forkmean} leads to
\[ y_{n+1}= \frac{n-3+\al}{n-\al } y_n + \frac{2- \al}{n-\al } x_n, \quad n \ge 3.\]
Solving this recurrence relation gives us
$y_4=\frac{\al(5-3\al)}{2(3-\al)(3-2\al)}$ and
for $n\ge 5$
\begin{align*}
y_n
&=y_3 \prod_{i=3}^{n-1}\frac{i-3+\al}{i-\al }   + \sum_{i=3}^{n-1} \frac{2- \al}{i-\al } x_i \prod_{j=i+1}^{n-1}\frac{j-3+\al}{j-\al }\\
& =\frac{1}{2}  \prod_{i=3}^{n-1}\frac{i-3+\al}{i-\al } +  \frac{(2- \al) \al }{2(3-2\al)} \sum_{i=3}^{n-1} \prod_{j=i+1}^{n-1}\frac{j-3+\al}{j-\al } \times \frac{1}{i-\al}\times \prod_{j=3}^{i-1} \frac{j-2+\al}{j-\al} \\
& =\frac{1}{2} \prod_{i=3}^{n-1}\frac{i-3+\al}{i-\al } + \frac{(2- \al) \al}{2(3-2\al)} \sum_{i=3}^{n-1}  \frac{1}{3-\al}\prod_{j=4}^{n-1}\frac{j-3+\al}{j-\al } \\
& =\frac{1}{2}  \prod_{i=3}^{n-1}\frac{i-3+\al}{i-\al } + \frac{(2- \al) \al}{2(3-2\al)} \frac{(n-3)}{(3-\al)} \prod_{j=4}^{n-1}\frac{j-3+\al}{j-\al }\\
&=
\frac{\al(2n-3+\al-n\al) }{2(3-2\al)(3-\al)}  \prod_{j=4}^{n-1}\frac{j-3+\al}{j-\al }.
\end{align*}

By  Lemma \ref{Lemma2},
\begin{align*}
y_n  &= \frac{\al(2n-3+\al-n\al) }{2(3-2\al)(3-\al)} \frac{\Gamma(4-\al) }{\Gamma(1+\al)} n^{-3+2\al} \left(1+o(1) \right)\\
&= \frac{\al(2-\al) }{2(3-2\al)} \frac{\Gamma(3-\al) }{\Gamma(1+\al)} n^{-2+2\al} \left(1+o(1) \right),
\end{align*}
as $n \to \infty$.
This completes the proof of part (ii) and hence the theorem. 
\end{proof}

\text{In the remainder of this section we present the proof of the second theorem.}

\begin{proof}[Proof of Theorem~\ref{thm:secondm}]


\tw{
Since the theorem clearly holds for $\al=1$, we shall assume that $\al\in [0,1)$ in the remainder of the proof. Furthermore, we will use the same $x_n$ and $y_n$ as defined in Theorem~\ref{thm:firstm}, and the fact that $x_n=\OO(n^{-2(1-\al)})$ and $y_n=\OO(n^{-2(1-\al)})$ as $n\to \infty$}.
We start with the proof of Part (i). To this end, we let 
\begin{equation}
\label{CM2}
z_n=\Ebold[C_n^2]  -\frac{(1-\al)^2}{(3-2\al)^2} \ n^2 - \frac{2(1-\al)(1+ 2\al -2\al^2)}{(5-4\al)(3-2\al)^2} \ n + \frac{\al (8-17\al+8\al^2)}{4(5-4\al)(3-2\al)^2}, ~n\ge 3.
\end{equation} 
Since $\Ebold[C_3^2]=1$, we get 
$ z_3 = \frac{88 \al^3 - 213 \al^2 + 152 \al - 24}{4(3-2\al)^2 (5-4\al)}.$
Next, substituting (\ref{CM2}) into (\ref{cherry2nd}) 
leads to
$$ (n-\al) z_{n+1} - (n-4+3\al) z_n = 2(1-\al)(n-1) x_n, \quad n \ge 3.  $$
Furthermore, using Theorem~\ref{thm:firstm}(i) we have
\begin{equation}
\label{eq:thm2:i}
\var(C_n)=\Ebold[C_n^2]-(\Ebold[C_n])^2= \frac{(1-\al)(2-\al)}{(5-4\al)(3-2\al)^2} \ n -  \frac{\al (1-\al)(2-\al)}{(5-4\al)(3-2\al)^2} + v_n -x_n^2,
\end{equation}
where
 \[v_n=z_n  - \frac{2(1-\al)}{3-2\al} n x_n - \frac{\al}{3-2\al} x_n = z_n - \frac{[2(1-\al) n +\al]}{3-2\al}x_n.\]
Then,  for $n\geq 3$, we have
\begin{align*}
(n-\al) v_{n+1} 
&= (n-\al)z_{n+1}  - \frac{[2(1-\al) (n+1) +\al]}{3-2\al} (n-\al)x_{n+1}
= (n-4+3\al)v_n  -\frac{2(1-\al)}{3-2\al} x_n
\end{align*}
and hence also
\begin{equation}
\label{eq:thm2:rec:v}
  v_{n+1} =\frac{n-4+3\al}{n-\al} v_n  -\frac{2(1-\al)}{(3-2\al)} \frac{x_n}{(n-\al)}.
  \end{equation}
\tw{
Now consider $f_n = \frac{n-4+3\al}{n-\al}$ and  $g_n  = -\frac{2(1-\al)x_n}{(3-2\al)(n-\al)}$ for $n\ge 3$, and let $a = 4-4\al$ and  $b= 3-2\al$. Then by Lemma~\ref{Lemma2}, it follows that there exists a constant $K_1$ such that
 \begin{align*}
    \prod_{i=\ell}^{n-1} |f_i|
= \prod_{i=\ell}^{n-1} \frac{i-4+3\al}{i-\al}
\le  K_1\left(\frac{n}{\ell} \right)^{-4+4\al}
=K_1\left(\frac{n}{\ell} \right)^{-a}
\quad\text{for all}\quad 3\le \ell  \le n-1,
\end{align*}
and by Theorem~\ref{thm:firstm}  there exists a constant $K_2$ such that
$$
|g_n|=\frac{2(1-\al)x_n}{(3-2\al)(n-\al)}
< \frac{4(1-\al)}{(3-2\al)}\frac{x_n}{n}
\le K_2\,n^{-3+2\al}
=K_2 n^{-b}
\quad\text{for}\quad n\ge 3. 
$$
Since $a-b+1=2(1-\al)>0$ for $\al\in [0,1)$, 
an application of Lemma \ref{Lemma1}
on the recursion~\eqref{eq:thm2:rec:v}
with the above $f_n$, $g_n$, $a$, $b$, 
and $c=\max\{K_1,K_2\}$ leads to  $v_n = \OO(n^{-2+2\al})$, and hence also  $v_n-x_n^2=\OO(n^{-2+2\al})$.
This, together with~\eqref{eq:thm2:i}, completes the proof of (i).\\
}


To prove Part (ii), 
we consider
\begin{eqnarray*}
t_n=\Ebold[A_nC_n] - \frac{(1-\al)^2}{2(3-2\al)^2} \ n^2+\frac{(1-\al)(4-25\al+16\al^2)}{4(5-4\al)(3-2\al)^2} n  +\frac{\al (8- 17\al +8\al^2 )}{4(5-4\al)(3-2\al)^2} \label{Eq:E[AnCn]}
\end{eqnarray*}
for $n\ge 3$. Combining   (\ref{cov}) and \eqref{Eq:E[AnCn]} leads to 
\begin{equation*}
(n-\al)t_{n+1}-(n-5+3\al)t_n = (2-\al) z_n + (1-\al)(n-1)y_n, \quad n\geq 3.
\end{equation*}
Since $\cov(A_n, C_n) =\Ebold[A_nC_n] -\Ebold[A_n] \Ebold[C_n]$, by \eqref{Exp:mu}, \eqref{Exp:nu} and \eqref{Eq:E[AnCn]} we have
\begin{equation}
\label{eq:thm2:ii}
\cov(A_n, C_n) = \frac{-(1-\al)(2-\al)(1-2\al)}{2(5-4\al)(3-2\al)^2} \ n -\frac{ \al (1-\al)(2-\al)}{(5-4\al)(3-2\al)^2}  + w_n -x_n y_n,
\end{equation}
where $w_n = t_n  -\dfrac{[(1-\al)n +\al]x_n + [2(1-\al)n +\al]y_n}{2(3-2\al)}$. 
Using straightforward but tedious algebraic simplification steps, we can show that 
\begin{align*}
(n-\al) w_{n+1}  -  (n-5+3\al)w_n &= 
(2-\al)v_n - \frac{1-\al}{3-2\al} x_n
\end{align*}
holds for $n\geq 3$ and hence 
\begin{equation}\label{Rec:w}
 w_{n+1}  =\frac{ n-5+3\al}{n-\al} w_n + (2-\al) \frac{v_n}{n-\al} - \frac{(1-\al)}{(3-2\al)} \frac{x_n}{(n-\al)}.
\end{equation}
Similarly to the proof of Part (i), applying Lemma \ref{Lemma1} on the recursion~\eqref{Rec:w} with $f_n = \frac{n-5+3\al}{n-\al} $,  $g_n  = \ (2-\al) \frac{v_n}{n-\al} - \frac{(1-\al)}{(3-2\al)} \frac{x_n}{(n-\al)} $,  $a = 5-4\al$ and $b= 3-2\al$,  we get  $w_n = \OO(n^{-2+2\al})$, and hence $w_n-x_ny_n = \OO(n^{-2+2\al})$. 
This proves part (ii) in view of~\eqref{eq:thm2:ii}.\\


To prove (iii), 
we consider
\[s_n= \Ebold[A_n^2]  - \frac{(1-\al)^2}{4(3-2\al)^2} n^2 - \frac{2(1-\al)(1+ 2\al +2\al^2)}{(5-4\al)(3-2\al)^2} n -   \frac{\al( 5-3\al+\al^2)}{4(3-2\al) (5-4\al)(7-4 \al)}  \]
for $n\ge 3$. 
Let $u_n = s_n -  \dfrac{[(1-\al)n + \al] y_n}{3-2\al}$. Then,
by straightforward  simplification steps we have
\begin{equation}
\label{eq:thm2:rec:u}
 u_{n+1} =  \frac{n-6+3\al}{n-\al} u_n  +  \frac{\al(2-\al)}{(3-2\al)} \frac{y_n}{(n-\al)} +\frac{(2-\al)^2}{(3-2\al)}  \frac{x_n}{(n-\al)}, \quad n\geq 3.
 \end{equation}
Furthermore, by $\var(A_n)=\Ebold[A_n^2]-(\Ebold[A_n])^2$ and Theorem~\ref{thm:firstm}(ii) we have
 \begin{eqnarray*} \var(A_n) 
  &=& \frac{(1-\al)(69-135\al+96\al^2-24\al^3)}{4(3-2\al)^2(5-4\al)(7-4\al)} \ n  
 + \frac{3\al(1-\al)(1-2\al)(5-3\al)}{4(3-2\al)^2(5-4\al)(7-4\al)} 
 + {\color{blue} u_n} -y_n^2.
\end{eqnarray*}
Similarly to the proof of Part (i),  applying Lemma \ref{Lemma1} on the recursion~\eqref{eq:thm2:rec:u} with $f_n = \frac{n-6+3\al}{n-\al}$,  $g_n  =\left[\frac{\al(2-\al)}{(3-2\al)} \frac{y_n}{(n-\al)} +\frac{(2-\al)^2}{(3-2\al)}  \frac{x_n}{(n-\al)}\right]$,  $a = 6-4\al$ and $b= 3-2\al$,  we get  $u_n = \OO(n^{-2+2\al})$ and hence also $u_n-y_n^2=\OO(n^{-2+2\al})$, which completes the proof of (iii) and hence the theorem.
\end{proof}



\section*{Acknowledgements} 
The work of Gursharn Kaur was supported by NUS Research Grant R-155-000-198-114, and that of Kwok Pui Choi by Singapore Ministry of Education Academic Research Fund
R-155-000-188-114. We thank Chris Greenman and Ariadne Thompson for stimulating discussions on the Ford model.

\bibliographystyle{plain}      
\bibliography{0_StatPhyloTree.bib}

\end{document}